# ON THE FALSE DISCOVERY RATE AND AN ASYMPTOTICALLY OPTIMAL REJECTION CURVE[1]

By Helmut Finner, Thorsten Dickhaus and Markus Roters

*German Diabetes Center, German Diabetes Center and Omnicare Clinical Research*

In this paper we introduce and investigate a new rejection curve for asymptotic control of the false discovery rate (FDR) in multiple hypotheses testing problems. We first give a heuristic motivation for this new curve and propose some procedures related to it. Then we introduce a set of possible assumptions and give a unifying short proof of FDR control for procedures based on Simes' critical values, whereby certain types of dependency are allowed. This methodology of proof is then applied to other fixed rejection curves including the proposed new curve. Among others, we investigate the problem of finding least favorable parameter configurations such that the FDR becomes largest. We then derive a series of results concerning asymptotic FDR control for procedures based on the new curve and discuss several example procedures in more detail. A main result will be an asymptotic optimality statement for various procedures based on the new curve in the class of fixed rejection curves. Finally, we briefly discuss strict FDR control for a finite number of hypotheses.

**1. Introduction.** The false discovery rate (FDR) has become one of the main research objects in multiple decision theory during the last decade. One reason for the increasing interest in FDR controlling procedures is the need for a suitable error rate controlling criterion if the number of hypotheses is very large. It is widely accepted that strong control of the familywise error rate (FWER) is far too restrictive for large systems of hypotheses. Beginning with the pioneering paper [2] on a linear step-up FDR controlling procedure based on Simes' test (cf. [22]) for independent $p$-values, meanwhile a series

Received January 2007; revised July 2007.
[1]Supported by the Deutsche Forschungsgemeinschaft.

of alternative procedures not only for independent but also for dependent test statistics have been introduced and theoretically investigated. For the latest developments we refer to [1, 5, 6, 7, 8, 16, 18, 20].

In general, the original linear step-up (LSU) procedure does not exhaust the prespecified FDR level, which gives hope for improvement with respect to power. In this paper we are concerned with the problem of finding an in some sense optimal rejection curve when the number of hypotheses tends to infinity.

The paper is organized as follows: In Section 2 we introduce some notation around the FDR and a basic model which will be considered throughout the paper. Moreover, we briefly describe the interaction of rejection curves, critical value curves and the empirical cumulative distribution function (ecdf) of $p$-values. In Section 3 we give a heuristic derivation based on asymptotics of the new rejection curve denoted by $f_\alpha$ and introduce some possible procedures based on this curve. Section 4 is of more general character and deals with assumptions, methods of proof, least favorable configurations (LFCs) for FDR control and upper FDR bounds. We first introduce a series of possible assumptions and give a unifying short proof of FDR control for procedures based on Simes' critical values which, among others, summarizes the derivations in [4, 19] for dependent $p$-values (or test statistics) in an elegant way. Then we investigate the problem of finding LFCs for the FDR, that is, we look for distributional settings where the FDR becomes largest and derive upper bounds for the FDR. Section 5 is concerned with asymptotic FDR control for procedures based on $f_\alpha$ under the assumption of independent and uniformly distributed $p$-values on the unit interval with respect to the set of true null hypotheses. Moreover, we formalize the asymptotic optimality of $f_\alpha$. Some concluding remarks, including a brief discussion on some properties of procedures related to $f_\alpha$ if a fixed number of hypotheses is at hand, are given in Section 6.

**2. Notation and preliminaries.** Throughout the paper, we use the following notation: Let $(\mathcal{X}, \mathcal{A}, \{P_\vartheta : \vartheta \in \Theta\})$ denote a statistical experiment and let $(H_n)_{n \in \mathbb{N}}$ be a sequence of null hypotheses with $\varnothing \neq H_n \subset \Theta$. The corresponding alternatives are given by $K_n = \Theta \setminus H_n$. Let $(p_n)_{n \in \mathbb{N}}$ denote a sequence of $p$-values with $p_n : (\mathcal{X}, \mathcal{A}) \longrightarrow ([0,1], \mathcal{B})$, where $\mathcal{B}$ denotes the Borel-$\sigma$-field over $[0,1]$. Denote the set of positive integers by $\mathbb{N}$ and let $\mathbb{N}_n = \{1, \ldots, n\}$, $I_0 = I_0(\vartheta) = \{i \in \mathbb{N} : \vartheta \in H_i\}$, $I_1 = I_1(\vartheta) = \mathbb{N} \setminus I_0 = \{i \in \mathbb{N} : \vartheta \notin H_i\}$ and $I_{n,j} = I_{n,j}(\vartheta) = I_j \cap \mathbb{N}_n$, $j = 0, 1$. As usual, let a $p$-value $p_i$ for testing $H_i$ satisfy $0 < P_\vartheta(p_i \leq x) \leq x$ for all $\vartheta \in H_i$, $i \in \mathbb{N}$ and $x \in (0,1]$. We assume that the measurable space $(\mathcal{X}, \mathcal{A})$ is large enough to accommodate probability measures $P_{I_0}$, $I_0 \subseteq \mathbb{N}$, under which all $p$-values $p_i, i \in I_0$, are i.i.d. uniformly distributed on $[0,1]$, and all $p_i$, $i \in I_1$, follow a Dirac distribution with point mass 1 in 0. We refer to $P_{I_0}$ as a *Dirac-uniform configuration*. We also assume that for every $\vartheta \in \Theta$ and $i \in I_0(\vartheta)$ there is a probability measure $P_{\vartheta i}$ defined on $(\mathcal{X}, \mathcal{A})$ for which the sequence $(p_n)_{n \in \mathbb{N}}$ has the same distribution under $P_{\vartheta i}$ as the sequence $(p_n^i)_{n \in \mathbb{N}}$ under $P_\vartheta$, the only difference between these two sequences of $p$-values being that $p_i^i \equiv 0$. This is a technical assumption which will be used in Section 5 for the determination of upper bounds for the FDR. Notice that the $P_{I_0}$'s and $P_{\vartheta i}$'s need not be contained in $\{P_\vartheta : \vartheta \in \Theta\}$. Under each $P_{I_0}$, the *Extended Glivenko–Cantelli theorem* (cf. [21], page 105) applies for the ecdf $\hat{F}_n$ (say) of the $p$-values, that is,

$$(2.1) \qquad \lim_{n \to \infty} \sup_{t \in [0,1]} \left| \hat{F}_n(t) - \left( \frac{n_1(n)}{n} + \frac{n_0(n)}{n} t \right) \right| = 0 \quad [P_{I_0}],$$

where $n_j = n_j(n) = |I_{n,j}|$, $j = 0, 1$. For a nonempty set $I_0 \subseteq \mathbb{N}$ we denote by $I_0'$ the set $I_0 \setminus \{\min I_0\}$ in the sequel.

For notational convenience we introduce $\hat{F}_{n,j}$, $j = 0, 1$, as the ecdfs of the $p$-values corresponding to the true ($j = 0$) and the false ($j = 1$) hypotheses, respectively. Now let $\mathcal{H}_n = \{H_i : i \in \mathbb{N}_n\}$ and let $\varphi_{(n)} = (\varphi_i : i \in \mathbb{N}_n)$ denote a non-randomized multiple test procedure for $\mathcal{H}_n$. For a fixed $\vartheta \in \Theta$ and a given test $\varphi_{(n)}$ we define the number of false rejections by $V_n = |\{i \in \mathbb{N}_n : \varphi_i = 1 \text{ and } \vartheta \in H_i\}|$ and the number of all rejections by $R_n = |\{i \in \mathbb{N}_n : \varphi_i = 1\}|$. Then the actual (expected) FDR of a multiple test $\varphi_{(n)}$ given $\vartheta \in \Theta$ is defined by

$$\mathrm{FDR}_\vartheta(\varphi_{(n)}) = \mathbb{E}_\vartheta \left[ \frac{V_n}{R_n \vee 1} \right].$$

The ratio $V_n / [R_n \vee 1]$ is the false discovery proportion (FDP). A multiple test $\varphi_{(n)}$ is said to control the FDR at level $\alpha \in (0, 1)$, if

$$\sup_{\vartheta \in \Theta} \mathrm{FDR}_\vartheta(\varphi_{(n)}) \leq \alpha.$$

The original LSU procedure for independent $p$-values ($\varphi_{(n)}^{\mathrm{LSU}}$ for short) rejects $H_i$ iff $p_i \leq m\alpha/n$, where $m = \max\{i \in \mathbb{N}_n : p_{i:n} \leq \alpha_{i:n}^{\mathrm{LSU}}\}$ with $\alpha_{i:n}^{\mathrm{LSU}} = i\alpha/n$ for $i = 1, \ldots, n$ and $p_{1:n} \leq \cdots \leq p_{n:n}$ denoting the ordered $p$-values. This procedure can be rewritten in terms of the ecdf $\hat{F}_n$ of the $p$-values. Let $t(\varphi_{(n)}^{\mathrm{LSU}}) = \sup\{t \in [0,1] : \hat{F}_n(t) \geq t/\alpha\}$. Then $\varphi_{(n)}^{\mathrm{LSU}}$ rejects $H_i$ iff $p_i \leq t(\varphi_{(n)}^{\mathrm{LSU}})$. The rejection curve $r_\alpha(t) = t/\alpha$ will be called the Simes line. Notice that $\alpha_{i:n}^{\mathrm{LSU}} = r_\alpha^{-1}(i/n)$. Now let $\vartheta \in \Theta$ and suppose that the $p_i$'s, $i \in I_{n,0}(\vartheta)$, are i.i.d. uniformly distributed on $[0, 1]$. Then one of the most interesting results for the LSU procedure is that

$$(2.2) \qquad \mathrm{FDR}_\vartheta(\varphi_{(n)}^{\mathrm{LSU}}) = \frac{n_0}{n} \alpha.$$

Different proofs of this fact can be found in [4, 11, 19, 23]. In [2], it was only shown that the FDR is always bounded by $n_0 \alpha/n$ for the LSU procedure. However, the fact that the FDR is bounded by $n_0 \alpha/n$ raised hope that improvements of the LSU procedure should be possible. In [13], an improvement taking a Bayes risk approach has been developed. Another idea is to estimate the number $n_0$ of true null hypotheses and to adjust the LSU procedure in order to exhaust the prespecified FDR level. This approach has been worked out in [23]. Multi-stage adaptive methods were recently proposed in [3]. The introduction in [14] gives a good overview about the development of several (adaptive) approaches aiming at improving the LSU procedure. In this paper we tackle the problem in a more direct way by deriving a new rejection curve which will serve as a basis for various stepwise multiple test procedures.

Many multiple test procedures can be described in terms of the ecdf $\hat{F}_n$ of the $p$-values and a rejection curve $r$. Let $\rho:[0,1] \to [0,1]$ be nondecreasing and continuous with $\rho(0) = 0$ and positive values on $(0,1]$. Define critical values $\alpha_{i:n} = \rho(i/n) \in (0,1]$ for $i = 1, \ldots, n$. We call $\rho$ a critical value function. Moreover, $r$ defined by $r(x) = \inf\{u \in [0,1] : \rho(u) = x\}$ for $x \in [0,1]$ ($\inf \varnothing = \infty$), will be called a rejection curve. For illustrative purposes, a plot of $\hat{F}_n$ together with the rejection curve $r$ is most useful in order to demonstrate the decision procedure. Note that we have the following relationship between the ecdf $\hat{F}_n$ of $n$ distinct $p$-values $p_1, \ldots, p_n$, the ordered $p$-values, the critical values $\alpha_{i:n} = \rho(i/n)$ and the rejection curve $r$:

$$\hat{F}_n(p_{i:n}) \geq r(p_{i:n}) \quad \text{if and only if} \quad p_{i:n} \leq \alpha_{i:n}.$$

A point $t = p_{i:n}$ satisfying $\hat{F}_n(p_{i:n}) \geq r(p_{i:n})$ and $\hat{F}_n(p_{i+1:n}) < r(p_{i+1:n})$ is called a crossing point (CP) between $\hat{F}_n$ and $r$. In this paper, we consider multiple test procedures where one of the CPs is chosen as a threshold $t^*$ in order to reject all $H_i$ with $p_i \leq t^*$. Other thresholding rules are extensively studied in [1]. It is immediately clear that the proportion $(R_n - V_n)/(n_1 \vee 1)$ of rejected false null hypotheses with respect to all false null hypotheses is nondecreasing in the threshold $t^*$. Therefore, we look for procedures which maximize $t^*$ for any given set of $p$-values subject to FDR control. Often, $\mathbb{E}_\vartheta[(R_n - V_n)/(n_1 \vee 1)]$ is defined as power. Loosely speaking, larger thresholds lead to larger power.

**3. Heuristic derivation of a new rejection curve.** In order to derive a new rejection curve we consider Dirac-uniform configurations $P_{I_0}$, $I_0 \subseteq \mathbb{N}$. This policy is motivated by the fact that Dirac-uniform configurations often provide upper bounds for the FDR, see Section 4. Let $\zeta_n = n_0(n)/n$ denote the proportion of true hypotheses among the first $n$ hypotheses. We refer to this situation as the Dirac-uniform finite model, $\mathrm{DU}_n(\zeta_n)$. Now suppose that

$$\lim_{n\to\infty} \zeta_n = \zeta \in [0,1].$$

Then (2.1) implies that, for $n$ tending to infinity, the ecdf $\hat{F}_n$ of the observed $p$-values converges to

$$F_\infty(t|\zeta) = (1-\zeta) + \zeta t \qquad \text{for all } t \in [0,1] \quad [P_{I_0}].$$

This situation will be called the Dirac-uniform asymptotic model, $\mathrm{DU}_\infty(\zeta)$ for short. Now suppose we choose some $t \in (0,1]$ and reject all hypotheses $H_i$ with $p_i \leq t$. Then, in the $\mathrm{DU}_\infty(\zeta)$ model, an expected proportion of $t\zeta$ of the true hypotheses will be rejected, that is, the false discovery rate $\mathrm{FDR}_\zeta(t)$ (say) is given by

$$\mathrm{FDR}_\zeta(t) = \frac{t\zeta}{(1-\zeta)+t\zeta}.$$

Now we determine a $t_\zeta \in [0,1]$ such that $\mathrm{FDR}_\zeta(t_\zeta)$ is exactly $\alpha$ for some prespecified $\alpha \in (0,1)$. Obviously, this can only work for $\zeta \in [\alpha, 1)$. In this case we obtain

$$t_\zeta = \frac{\alpha(1-\zeta)}{\zeta(1-\alpha)}.$$

For $\zeta \in [0,\alpha)$ one may set $t_\zeta = 1$, which implies that all hypotheses are rejected and $\mathrm{FDR}_\zeta(1) = \zeta < \alpha$.

We now look for a strictly increasing rejection curve $f_\alpha : [0,1] \to [0,1]$ which does not depend on $\zeta$ but tells us which of the hypotheses can be rejected such that the limiting FDR is exactly $\alpha$ for $\zeta \in [\alpha, 1)$. The function $f_\alpha$ can be found by requiring $f_\alpha(t_\zeta) = F_\infty(t_\zeta)$ for all $\zeta \in [\alpha, 1)$, which is equivalent to

$$f_\alpha\left(\frac{\alpha(1-\zeta)}{\zeta(1-\alpha)}\right) = \frac{1-\zeta}{1-\alpha}.$$

Obviously,

$$f_\alpha(t) = \frac{t}{t(1-\alpha)+\alpha}, \qquad t \in [0,1],$$

fulfills this requirement. In the $\mathrm{DU}_\infty(\zeta)$ model, the limiting ecdf $F_\infty(\cdot|\zeta)$ crosses $f_\alpha$ in the point $(t_\zeta, y_\zeta)$ with $t_\zeta = \min\{\alpha(1-\zeta)/[\zeta(1-\alpha)], 1\}$ and $y_\zeta = \min\{(1-\zeta)/(1-\alpha), 1\}$ for the first time, and it should be noted that this is also true for $\zeta = 1$ if we define $(t_1, y_1) = (0,0)$.

REMARK 3.1. Comparing the Simes line $r_\alpha(t) = t/\alpha$ and the new rejection curve $f_\alpha(t)$, we obviously have $r_\alpha(t) > f_\alpha(t)$ for $t > 0$ and the derivative in $t = 0$ is $1/\alpha$ for both curves. Moreover, notice that $f_\alpha$ obeys the symmetry property

$$f_\alpha^{-1}(t) = \frac{\alpha t}{1 - (1-\alpha)t} = 1 - f_\alpha(1 - t) \qquad \text{for all } t \in [0, 1],$$

where $f_\alpha^{-1}$ denotes the inverse of $f_\alpha$. Clearly, $f_\alpha^{-1}$ is a critical value function.

The new rejection curve $f_\alpha$ will be called the asymptotically optimal rejection curve (AORC). The question is how to implement $f_\alpha$ not only in the Dirac-uniform models but also in more general models into a multiple test procedure which controls the FDR level $\alpha$ strictly or at least asymptotically. The critical values induced by $f_\alpha$ are given by

$$(3.1) \quad \alpha_{i:n} = f_\alpha^{-1}\left(\frac{i}{n}\right) = \frac{(i/n)\alpha}{1 - (i/n)(1 - \alpha)} = \frac{i\alpha}{n - i(1 - \alpha)}, \qquad i = 1, \ldots, n.$$

Remember that distinct $p$-values satisfy

$$\hat{F}_n(p_{i:n}) \geq f_\alpha(p_{i:n}) \quad \text{if and only if} \quad p_{i:n} \leq \alpha_{i:n}.$$

It is tempting to use $\alpha_{1:n} \leq \cdots \leq \alpha_{n:n}$ in a step-up (SU) procedure for testing $n$ hypotheses. Unfortunately, $\alpha_{n:n} = 1$, so that this procedure always rejects all hypotheses. This pitfall is due to the fact that $f_\alpha(1) = \hat{F}_n(1) = 1$. Therefore, we need some adjustment with respect to $f_\alpha$ or the SU procedure. In the remainder of this section, we consider some candidates for asymptotic FDR control avoiding the aforementioned pitfall.

EXAMPLE 3.1 (Step-up-down procedures). An interesting class of procedures based on critical values $0 < \alpha_{1:n} \leq \cdots \leq \alpha_{n:n} \leq 1$ are step-up-down (SUD) procedures introduced in [24] and studied in [19] in terms of FDR control. For $\lambda_n \in \mathbb{N}_n$ an SUD procedure $\varphi_{n,\lambda_n}^{\text{SUD}} = (\varphi_1, \ldots, \varphi_n)$ of order $\lambda_n$ is defined as follows: If $p_{\lambda_n:n} \leq \alpha_{\lambda_n:n}$, set $m_n = \max\{j \in \{\lambda_n, \ldots, n\} : p_{i:n} \leq \alpha_{i:n}$ for all $i \in \{\lambda_n, \ldots, j\}\}$, whereas for $p_{\lambda_n:n} > \alpha_{\lambda_n:n}$, put $m_n = \sup\{j \in \{1, \ldots, \lambda_n\} : p_{j:n} \leq \alpha_{j:n}\}$ ($\sup \varnothing = -\infty$). Define $\varphi_i = 1$ if $p_i \leq \alpha_{m_n:n}$ and $\varphi_i = 0$ otherwise ($\alpha_{-\infty:n} = -\infty$). Note that $\lambda_n = 1$ yields a step-down (SD) procedure and $\lambda_n = n$ yields an SU procedure. The order of an SUD procedure can be defined in terms of a fixed parameter $\lambda \in [0, 1]$ by setting $\lambda_n = \inf\{j \in \mathbb{N}_n : \alpha_{j:n} \geq \lambda\}$ ($\inf \varnothing = n$). Then $\lambda = 0$ ($\lambda = 1$) corresponds to an SD (SU) procedure. An SUD procedure of order $\lambda_n = \lambda_n(\lambda), \lambda \in [0, 1)$, based on $f_\alpha$ resolves the problems around the point $t = 1$ in an elegant way. It is obvious in view of Remark 3.1 that in the case of $\lambda \geq \alpha$ the new SUD procedure based on $f_\alpha$ rejects at least all hypotheses rejected by the LSU procedure, possibly more. Therefore, one cannot expect that the FDR level is controlled in the finite case. However, it will be shown that the FDR is controlled asymptotically. Note that $\varphi_{n,\lambda_n}^{\text{SUD}}$ is component-wise nondecreasing in $\lambda$.

EXAMPLE 3.2 (Adjusted SU procedures based on $f_\alpha$). As noted before, an SU procedure based on $\alpha_{i:n} = f_\alpha^{-1}(i/n)$ cannot work. Therefore, some adjustment of $f_\alpha$ in an SU procedure is necessary. We first consider the case where the adjusted rejection curve $f_\alpha^{\text{adj}}$ satisfies that $f_\alpha^{\text{adj}}(x)/x$ is nonincreasing in $x$, an important property for the calculation of the FDR. One may specify some $\kappa \in (0,1)$ and define a new rejection curve

$$f_{\alpha,\kappa}^{\text{adj}}(x) = \begin{cases} f_\alpha(x), & 0 \leq x < \kappa, \\ h(x), & \kappa \leq x \leq 1, \end{cases}$$

such that $f_\alpha^{\text{adj}}(x)/x$ is nonincreasing in $x$, $f_\alpha^{\text{adj}}(x^*) = 1$ for some $x^* < 1$, and, $f_\alpha(\kappa) = h(\kappa)$. For example, one may choose $h(x) = ax + b$ for suitable values $a,b \geq 0$. We consider two possible choices of $h$ ($h_1$ and $h_2$ say) and refer to the resulting rejection curves as $f_{\alpha,\kappa}^{(i)}, i = 1,2$. Let

$$h_1(x) = f_\alpha'(\kappa)(x - \kappa) + f_\alpha(\kappa), \qquad x \in [\kappa, 1].$$

Then $h_1'(\kappa) = f_\alpha'(\kappa)$, $h_1(\kappa) = f_\alpha(\kappa)$ and $h_1(x^*) = 1$ for $x^* = \kappa(1-\alpha)(2-\kappa) + \alpha < 1$. The largest possible slope of $h$ is $a = f_\alpha(\kappa)/\kappa$. This leads to the second choice, that is, $h_2(x) = x f_\alpha(\kappa)/\kappa$. This is close to the truncated SU procedure in Example 3.3 below. Note that $h_2(x^*) = 1$ for $x^* = \kappa(1-\alpha) + \alpha$. For example, suppose that $\kappa = f_\alpha^{-1}(i/n)$ for some fixed $i \in \mathbb{N}_n$. Then the step-up critical values are given by

$$\gamma_{j:n} = \begin{cases} f_\alpha^{-1}(j/n), & 1 \leq j \leq i, \\ \dfrac{j}{n} \dfrac{\kappa}{f_\alpha(\kappa)}, & i+1 \leq j \leq n. \end{cases}$$

EXAMPLE 3.3 (Truncated SU procedures based on $f_\alpha$). Let $\kappa \in (0,1)$ be fixed and define

$$\rho_\alpha(x) = \begin{cases} f_\alpha^{-1}(x), & 0 \leq x \leq f_\alpha(\kappa), \\ \kappa, & f_\alpha(\kappa) < x \leq 1. \end{cases}$$

With $\gamma_{i:n} = \min\{f_\alpha^{-1}(i/n), \kappa\}$ we have $\gamma_{j:n} = \rho_\alpha(i/n)$ for $j = i, \ldots, n$. Hence, the truncated SU procedure is well defined in terms of $\rho_\alpha$. It is worth mentioning that this type of procedure differs substantially from the adjusted procedures discussed in Examples 3.1 and 3.2, because the monotonicity behavior of the ratio $\rho_\alpha(x)/x$ changes at $x = f_\alpha(\kappa)$, which makes FDR calculation much subtler.

The next section deals with methodologies for proving FDR control and deriving upper FDR bounds. The results will be applied in Section 5 in order to derive conditions for asymptotic FDR control for procedures based on the AORC. The aforementioned example procedures are then investigated at the end of Section 5 with respect to asymptotic FDR control and asymptotic optimality.

**4. A unifying proof of FDR control and upper FDR bounds.** Suppose that $R_n$ and $\varphi_{(n)}$, respectively, are defined in terms of $p$-values $p_1, \ldots, p_n$ and critical values $\alpha_{i:n} = \rho(i/n)$ for some critical value function $\rho$. Consider the following three sets of possible assumptions. The first two assumptions concern the structure of the test procedure (test assumptions):

(T1) $\forall i \in \mathbb{N}_n : p_i \leq \alpha_{1:n}$ implies $\varphi_i = 1$.
(T2) $\forall j \in \mathbb{N}_n : R_n = j$ implies $\forall i \in \mathbb{N}_n : [\varphi_i = 1 \Leftrightarrow p_i \leq \alpha_{j:n}]$.

The second set of assumptions concerns properties of distributions of $p$-values and $R_n$ (distributional assumptions):

(D1) $\forall \vartheta \in \Theta : \forall j \in \mathbb{N}_n : \forall i \in I_{n,0}(\vartheta) : P_\vartheta(R_n \geq j | p_i \leq t)$ is nonincreasing in $t \in (0, \alpha_{j:n}]$.
(D2) $\forall \vartheta \in \Theta : \forall j \in \mathbb{N}_n : \forall i \in I_{n,0}(\vartheta) : \forall t \in (0, \alpha_{j:n}] : P_\vartheta(R_n \geq j | p_i \leq t) \leq P_{\vartheta^i}(R_n \geq j)$.
(D3) $\forall \vartheta \in \Theta : \forall i \in I_{n,0}(\vartheta) : p_i \sim U([0,1])$.

Finally we have two possible independence assumptions:

(I1) $\forall \vartheta \in \Theta$: The $p_i$'s, $i \in I_{n,0}(\vartheta)$, are i.i.d.
(I2) $\forall \vartheta \in \Theta : (p_i : i \in I_{n,0})$ and $(p_i : i \in I_{n,1})$ are stochastically independent random vectors.

Furthermore, we often refer to step-up tests or step-up-down tests of order $\lambda_n$. The simple structure of SU tests often simplifies derivations concerning properties of these tests. If $\varphi_{(n)}$ is a step-up-down test, the properties of a step-up test remain valid in the step-up branch of such a procedure. For example, it is important to note (see [19], page 248) that given a step-up-down test of order $\lambda_n$, under (D3), (I1) and (I2), we get for all $\vartheta \in \Theta$ and all $i \in I_{n,0}(\vartheta)$

(4.1) $\quad \forall j = 1, \ldots, \lambda_n : \forall t \in (0, \alpha_{j:n}] : P_\vartheta(R_n \geq j | p_i \leq t) = P_{\vartheta^i}(R_n \geq j),$

(4.2) $\quad \forall j = 1, \ldots, \lambda_n : \forall t \in (0, \alpha_{j:n}] : P_\vartheta(R_n = j | p_i \leq t) = P_{\vartheta^i}(R_n = j).$

For $\lambda_n = n$, that is, for a step-up test, we even get

$$\forall j = 1, \ldots, n : \forall t \in (0, \alpha_{j:n}] : P_\vartheta(R_n \geq j | p_i \leq t) = P_{\vartheta^i}(R_n \geq j).$$

Assumptions (T1) and (T2) concern more general structures of test procedures. Step-up-down tests satisfy both (T1) and (T2). The monotonicity assumption in (D1) is somewhat weaker than the PRDS assumption (PRDS: positive regression dependency on subsets). More precisely, from the $P_\vartheta^{p_i}$-almost sure antitonicity of the factorized conditional probability $P_\vartheta(R_n \geq j | p_i = t)$ in $t \in [0, \alpha_{j:n}]$ we obtain the property formulated in (D1), where the equality in the condition is replaced by an inequality. This type of conclusion is indicated in [17] and can be proved in an easy way by using Wijsman's inequality, see [25]. Anyhow, (D1) is the decisive condition for dependent $p$-values in order to prove FDR control or to derive upper bounds for the FDR. Examples of distributions being PRDS are extensively studied in [4] and [19]. Anyhow, we have no example with dependent $p$-values yet where (D1) is easier to show than PRDS. However, under the independence assumptions (I1) and (I2), the important class of SUD tests satisfy (D1). Property (D2) will only be used under the independence assumptions (I1) and (I2) and is an important tool for deriving LFC results. In case of dependency, (D2) is often violated. Assumptions (D3) and (I1) concern the distribution of $p$-values under the corresponding null hypotheses.

The following theorem and its proof unify, simplify and slightly extend the results and the proofs given in [4] and [19], respectively.

THEOREM 4.1. *Let $\alpha \in (0,1)$ and let $\varphi_{(n)}$ be a multiple test procedure for $\mathcal{H}_n$ defined in terms of Simes' critical values $\alpha_{i:n} = i\alpha/n$, $i = 1, \ldots, n$. Let $\vartheta \in \Theta$ such that $n_0$ hypotheses are true and the remaining ones are false. If (T1), (T2) and (D1) are satisfied, then*

$$\mathrm{FDR}_\vartheta(\varphi_{(n)}) \leq \frac{n_0}{n}\alpha,$$

*with "=" if $\varphi_{(n)}$ is a step-up test and (D3), (I1) and (I2) are satisfied in addition.*

PROOF. Consider the following chain of (in)equalities:

$$\mathrm{FDR}_\vartheta(\varphi_{(n)}) = \sum_{i \in I_{n,0}(\vartheta)} \sum_{j=1}^n \frac{1}{j} P_\vartheta(R_n = j, \varphi_i = 1)$$

$$(4.3) \qquad = \sum_{i \in I_{n,0}(\vartheta)} \sum_{j=1}^n \frac{1}{j} P_\vartheta(p_i \leq \alpha_{j:n}) P_\vartheta(R_n = j | p_i \leq \alpha_{j:n})$$

$$(4.4) \qquad \leq \sum_{i \in I_{n,0}(\vartheta)} \sum_{j=1}^n \frac{\alpha_{j:n}}{j} P_\vartheta(R_n = j | p_i \leq \alpha_{j:n})$$

$$
\begin{aligned}
&\leq \sum_{i \in I_{n,0}(\vartheta)} \Bigg[ \alpha_{1:n} P_\vartheta(R_n \geq 1 | p_i \leq \alpha_{1:n}) \\
&\qquad\qquad + \sum_{j=2}^{n} \bigg[ \frac{\alpha_{j:n}}{j} - \frac{\alpha_{j-1:n}}{j-1} \bigg] P_\vartheta(R_n \geq j | p_i \leq \alpha_{j:n}) \Bigg]
\end{aligned}
\tag{4.5}
$$

$$
= \frac{n_0}{n} \alpha. \tag{4.6}
$$

Equation (4.3) holds under (T2), and "=" holds in (4.4) if (D3) holds. Inequality (4.5) holds under the assumption (D1) with "=" if $\varphi_{(n)}$ is a step-up test and (D3), (I1) and (I2) hold. Finally, (4.6) is a consequence of (T1). □

REMARK 4.1. The key step in the proof is (4.5), where $P_\vartheta(R_n \geq j | p_i \leq \alpha_{j-1:n})$ is replaced by $P_\vartheta(R_n \geq j | p_i \leq \alpha_{j:n})$ for $j = 2, \ldots, n$ according to assumption (D1). In case of dependency or in case of a non-step-up test, the difference between these quantities may sum up to a considerable amount, that is, the FDR may be much smaller than the upper bound $n_0 \alpha / n$ in such cases. For a detailed investigation of the latter phenomenon, see [9].

One of the main problems in order to ensure FDR control of a multiple test procedure is to find least favorable parameter configurations (LFCs), that is, parameter configurations under which the FDR for a given test procedure becomes largest. Obviously, LFCs are no issue for the LSU procedure if (D3), (I1) and (I2) hold true. To date it looks like SU procedures are easier to cope with than SD or SUD procedures. One reason for this is that Dirac-uniform configurations can often be viewed as least favorable for certain SU procedures. This fact is based on the following important result:

THEOREM 4.2 (Benjamini and Yekutieli (2001); cf. [4]). *An SU procedure with critical values $\alpha_{1:n} \leq \cdots \leq \alpha_{n:n}$ fulfilling (D3), (I1) and (I2) has the following properties:*

(1) *If the ratio $\alpha_{i:n}/i$ is increasing in $i$, as $(p_i : i \in I_{n,1})$ increases stochastically, the FDR decreases.*

(2) *If the ratio $\alpha_{i:n}/i$ is decreasing in $i$, as $(p_i : i \in I_{n,1})$ increases stochastically, the FDR increases.*

Hence, under the assumptions of Theorem 4.2, the Dirac-uniform configurations, where all $p$-values under alternatives are almost surely 0, can be viewed as LFCs if the ratio $\alpha_{i:n}/i$ is increasing in $i$. More precisely, on the parameter subspace, where exactly $n_0$ ($n_1$) hypotheses are true (false), the FDR becomes largest if the $p$-values under alternatives are almost surely 0.

Therefore, it suffices to consider all Dirac-uniform configurations in order to check whether the FDR is controlled at level $\alpha$.

Unfortunately, the method of proof given in [4] does not seem to work for SD and SUD procedures. However, we show below that Dirac-uniform configurations often provide upper bounds. To this end, we define $q(x) = \rho(x)/x$ for all $x \in (0, 1]$ and assume that $q(0) = \limsup_{x \downarrow 0} q(x) < \infty$. Moreover, we define the (continuous) function $\overline{q}$ by $\overline{q}(x) = \max_{0 \leq t \leq x} q(t)$, $x \in [0, 1]$. Hence, $\overline{q}$ is the upper isotonic envelope or, in other words, the least isotonic majorant of $q$. For the derivation of upper FDR bounds, we now introduce the following additional conditions:

(A1) If $(p_1, \ldots, p_n)$ is stochastically not greater under $\vartheta_1 \in \Theta$ than under $\vartheta_2 \in \Theta$, then $\varphi_{(n)}$ is stochastically not greater under $\vartheta_2 \in \Theta$ than under $\vartheta_1 \in \Theta$.
(A2) $\rho(x)/x$ is nondecreasing for $x \in (0, 1]$.

Note that $\alpha_{i:n}/i$ is nondecreasing in $i$ if (A2) holds. If $\rho$ is differentiable on $(0,1)$, (A2) is equivalent to $\rho'(x) \geq q(x)$ for $x \in (0,1)$. In what follows, $\overline{q}$ is essential in deriving upper bounds for the FDR. Note that $q \neq \overline{q}$ for the truncated step-up procedure introduced in Example 3.3. If $q \neq \overline{q}$, the bounds for the FDR based on $\overline{q}$ may not be that sharp.

THEOREM 4.3. *Let $\vartheta \in \Theta$ such that $n_0 \in \mathbb{N}_n$, hypotheses are true and the remaining ones are false. Let $i_0 = \min I_0$ [and $I_0' = I_0 \setminus \{i_0\}$ as defined below equation (2.1)]. If (T1)–(I2) are satisfied, then*

$$\text{(4.7)} \qquad \text{FDR}_\vartheta(\varphi_{(n)}) \leq \frac{n_0}{n} \sum_{j=1}^n \overline{q}(j/n) P_{\vartheta^{i_0}}(R_n/n = j/n)$$

$$\text{(4.8)} \qquad \qquad = \frac{n_0}{n} \mathbb{E}_{\vartheta^{i_0}} \overline{q}(R_n/n),$$

*with equality in (4.7) if $\varphi_{(n)}$ is a step-up test and (A2) holds in addition. If (T1)–(I2) and (A1) are fulfilled, then*

$$\text{(4.9)} \qquad \text{FDR}_\vartheta(\varphi_{(n)}) \leq \frac{n_0}{n} \mathbb{E}_{I_0'} \overline{q}(R_n/n).$$

PROOF. Let $b_j = P_\vartheta(R_n \geq j | p_{i_0} \leq \alpha_{j:n})$ and $\Delta \overline{q}(j/n) = \overline{q}(j/n) - \overline{q}((j-1)/n)$ for $j = 1, \ldots, n$. Then, proceeding as in the proof of Theorem 4.1 we get for fixed $\vartheta \in \Theta$ under (D1)–(D3), (I1) and (I2)

$$\text{FDR}_\vartheta(\varphi_{(n)}) = \frac{n_0}{n} \sum_{j=1}^n q(j/n) P_\vartheta(R_n = j | p_{i_0} \leq \alpha_{j:n})$$

$$\text{(4.10)} \qquad \qquad \leq \frac{n_0}{n} \sum_{j=1}^n \overline{q}(j/n) P_\vartheta(R_n = j | p_{i_0} \leq \alpha_{j:n})$$

$$
(4.11) \qquad \leq \frac{n_0}{n}\left[\overline{q}(1/n)b_1 + \sum_{j=2}^{n}\Delta\overline{q}(j/n)b_j\right]
$$

$$
(4.12) \qquad \leq \frac{n_0}{n}\left[\overline{q}(1/n)P_{\vartheta^{i_0}}(R_n \geq 1) + \sum_{j=2}^{n}\Delta\overline{q}(j/n)P_{\vartheta^{i_0}}(R_n \geq j)\right]
$$

$$
= \frac{n_0}{n}\sum_{j=1}^{n}\overline{q}(j/n)P_{\vartheta^{i_0}}(R_n/n = j/n),
$$

which proves (4.7). In view of $P_{\vartheta^{i_0}}(R_n > 0) = 1$ according to (T1), (4.8) follows immediately. If $\varphi_{(n)}$ is a step-up test, which implies (4.1) for $\lambda_n = n$, we have equality in (4.11) and (4.12), hence $q = \overline{q}$ yields equality in (4.10). Finally, in order to prove (4.9), we use the same argumentation as in the proof of Theorem 4.2 given in [4], that is, that stochastic increase in the distribution of the random vector $(p_1, \ldots, p_n)$ can be characterized by the increase of the expectation of all nondecreasing functions (in case the expectation exists). To this end, we note that obviously $R_n = |\{i \in \mathbb{N}_n : \varphi_i = 1\}|$ is a nondecreasing function of $\varphi_{(n)}$ and hence, due to (A1), is stochastically nonincreasing in $(p_1, \ldots, p_n)$. The isotonicity of $\overline{q}$ completes the proof. □

Inequality (4.9) will be a helpful tool in order to calculate upper FDR bounds and to prove FDR control, because it only makes use of the distribution of $R_n$ under Dirac-uniform configurations. Especially for SUD tests, this distribution can be handled analytically.

**5. Asymptotic FDR control for procedures based on the AORC.** This section deals with conditions for asymptotic FDR control for procedures based on the new rejection curve. A major result will be that the example procedures presented in Section 3 control the FDR asymptotically. Theorems 5.1 and 5.3 provide sufficient conditions for asymptotic FDR control. If the underlying procedure leads to a determinable proportion of rejected hypotheses, Theorems 5.2 and 5.4 even give explicit values for the resulting FDR. Moreover, the asymptotic optimality of $f_\alpha$ is formalized in Theorem 5.5.

THEOREM 5.1. *Suppose $\varphi_{(n)}$ is based on $\rho \leq f_\alpha^{-1}$ and that (T1)–(I2) and (A1) are fulfilled. If for all nonempty sets $I_0 \subseteq \mathbb{N}$ and all subsequences $(n_k)_{k \in \mathbb{N}} \subseteq \mathbb{N}$ with $\lim_{k \to \infty} \zeta_{n_k} = \zeta$ for some $\zeta \in [0,1]$, it holds that*

$$
(5.1) \qquad \limsup_{k \to \infty} \frac{R_{n_k}}{n_k} \leq f_\alpha(t_\zeta) \quad [P_{I_0'}],
$$

*then*

$$
(5.2) \qquad \limsup_{n \to \infty} \sup_{\vartheta \in \Theta} \mathrm{FDR}_\vartheta(\varphi_{(n)}) \leq \alpha.
$$

PROOF. Let, for notational convenience, $P_{m,n}$ refer to a Dirac-uniform configuration such that the first $m$ $p$-values are i.i.d. uniformly distributed and the remaining ones follow a Dirac distribution with point mass in 0, $0 \leq m \leq n$, $n \in \mathbb{N}$. Then we have from inequality (4.9)

$$\forall n \in \mathbb{N}: \sup_{\vartheta \in \Theta} \mathrm{FDR}_\vartheta(\varphi_{(n)}) \leq \max_{1 \leq n_0 \leq n} \frac{n_0}{n} \mathbb{E}_{n_0-1,n} \overline{q}(R_n/n).$$

Since for each $n \in \mathbb{N}$ the maximum in this inequality is attained at some value $n_0(n)$ (say), we get

$$\limsup_{n \to \infty} \sup_{\vartheta \in \Theta} \mathrm{FDR}_\vartheta(\varphi_{(n)}) \leq \limsup_{n \to \infty} \zeta_n \mathbb{E}_{n_0(n)-1,n} \overline{q}(R_n/n),$$

where $\zeta_n = n_0(n)/n$, $n \in \mathbb{N}$. We now may extract a subsequence $(n_k)_{k \in \mathbb{N}}$ of $\mathbb{N}$ with $\lim_{k \to \infty} \zeta_{n_k} = \zeta$ for some $\zeta \in [0,1]$ such that

$$\limsup_{n \to \infty} \zeta_n \mathbb{E}_{n_0(n)-1,n} \overline{q}(R_n/n) = \lim_{k \to \infty} \zeta_{n_k} \mathbb{E}_{n_0(n_k)-1,n_k} \overline{q}(R_{n_k}/n_k)$$

$$\leq \zeta \limsup_{k \to \infty} \mathbb{E}_{n_0(n_k)-1,n_k} \overline{q}^*(R_{n_k}/n_k),$$

where $\overline{q}^*$ denotes the $\overline{q}$-function corresponding to the critical value function $f_\alpha^{-1}$. Similarly as in [11], pages 1003–1004, we are able to select from $(n_k)_{k \in \mathbb{N}}$ a further subsequence (without loss of generality with the same name) and construct a global set $I \subseteq \mathbb{N}$ with the property $|I \cap I_{n_k}| = n_0(n_k)$ for all $k \in \mathbb{N}$. (At this point it should be noted that the definition of the sets $M_k$ at the bottom of page 1003 in [11] has a typo at its right end in that the term $k(n_k)$ has to be replaced by $n_k$.) Now we obtain from (5.1)

$$\zeta \limsup_{k \to \infty} \mathbb{E}_{n_0(n_k)-1,n_k} \overline{q}^*(R_{n_k}/n_k) = \zeta \limsup_{k \to \infty} \mathbb{E}_{I'} \overline{q}^*(R_{n_k}/n_k)$$

$$\leq \zeta \mathbb{E}_{I'} \overline{q}^*\left(\limsup_{k \to \infty} R_{n_k}/n_k\right)$$

$$= \zeta \overline{q}^*(f_\alpha(t_\zeta))$$

$$= \min\{\alpha, \zeta\} \leq \alpha,$$

hence the assertion of the theorem, that is, (5.2) follows. □

If we sharpen assumption (5.1), we can even give explicit values for the FDR.

THEOREM 5.2. *Let $\vartheta \in \Theta$, $\varphi_{(n)}$ be based on $\rho \leq f_\alpha^{-1}$ and assume (T2), (D3), (I1) and*

(5.3) $$\lim_{n \to \infty} \zeta_n = \zeta \in [0,1].$$

*If $\lim_{n \to \infty} R_n/n = r^*[P_\vartheta]$ for some $r^* \in (0, f_\alpha(t_\zeta)]$, then it holds*

(5.4) $$\lim_{n \to \infty} \mathrm{FDR}_\vartheta(\varphi_{(n)}) = \zeta \rho(r^*)/r^* = \zeta q(r^*) \leq \min\{\alpha, \zeta\}.$$

PROOF. From (T2) and for $n_0, n \in \mathbb{N}$ we get the representation

$$V_n = n_0 \hat{F}_{n,0}(\rho(R_n/n))\mathbf{1}_{\{R_n > 0\}}.$$

From this we obtain the inequality chain

$$|V_n/n - \zeta_n \rho(R_n/n)| \leq \zeta_n |\hat{F}_{n,0}(\rho(R_n/n)) - \rho(R_n/n)| \leq \sup_{t \in [0,1]} |\hat{F}_{n,0}(t) - t|.$$

Hence, using the Glivenko–Cantelli property (2.1) together with the remaining assumptions of the theorem and the continuity of $\rho$, we finally see that $V_n/n$ converges $P_\vartheta$-almost surely to $\zeta \rho(r^*)$. Thus, due to $r^* > 0$, we have $\lim_{n \to \infty} \mathbb{E}_\vartheta[V_n/(R_n \vee 1)] = \zeta \rho(r^*)/r^*$. The right-hand side inequality in (5.4) is obtained by noting that $\zeta f_\alpha^{-1}(t)/t$ is increasing in $t \in (0, f_\alpha(t_\zeta)]$ to $\zeta t_\zeta/f_\alpha(t_\zeta) = \min\{\alpha, \zeta\}$ at $t = f_\alpha(t_\zeta)$. □

The remaining case $r^* = 0$ is treated in the following two theorems.

THEOREM 5.3. *Let $\vartheta \in \Theta$, $\varphi_{(n)}$ be based on $\rho \leq f_\alpha^{-1}$ and assume (T1)–(I2), (A1), (5.3) and*

(5.5) $$\forall \varepsilon > 0 : \liminf_{n \to \infty} \inf_{\varepsilon \leq t \leq 1} (t - \hat{F}_n(\rho(t))) > 0 \quad [P_\vartheta].$$

*Then it holds that*

(5.6) $$\limsup_{n \to \infty} \text{FDR}_\vartheta(\varphi_{(n)}) \leq \zeta \limsup_{x \downarrow 0} q(x) = \zeta q(0) = \zeta \overline{q}(0) \leq \zeta \alpha.$$

PROOF. To avoid triviality, we assume $I_0(\vartheta) \neq \varnothing$. Then, from (4.7) and (4.8) we have that

(5.7) $$\limsup_{n \to \infty} \text{FDR}_\vartheta(\varphi_{(n)}) \leq \zeta \limsup_{n \to \infty} \mathbb{E}_{\vartheta^{i_0}} \overline{q}(R_n/n).$$

Since due to (T1) and (T2) we have for all $n \in \mathbb{N}$

$$\hat{F}_n(\rho(R_n/n)) = R_n/n,$$

(5.5) implies that for every fixed $\varepsilon > 0$ we obtain $\limsup_{n \to \infty} R_n/n \leq \varepsilon$ $P_\vartheta$-almost surely, that is, $\lim_{n \to \infty} R_n/n = 0$ $P_\vartheta$-almost surely. Now, since for all $n \in \mathbb{N}$ the maximum absolute difference on the unit interval of the ecdf $\hat{F}_n$ [corresponding to the sequence of $p$-values $(p_n)_{n \in \mathbb{N}}$] and the ecdf $\hat{F}_n^{i_0}$ [corresponding to the sequence of $p$-values $(p_n^{i_0})_{n \in \mathbb{N}}$ defined in Section 2] is at most $1/n$, condition (5.5) also holds $P_{\vartheta^{i_0}}$-almost surely, which entails that $\lim_{n \to \infty} R_n/n = 0$ $P_{\vartheta^{i_0}}$-almost surely. Hence, due to the continuity of $\overline{q}$ we have $\lim_{n \to \infty} \mathbb{E}_{\vartheta^{i_0}} \overline{q}(R_n/n) = \overline{q}(0) = q(0) \leq \lim_{t \downarrow 0} f_\alpha^{-1}(t)/t = \alpha$. In view of inequality (5.7), this completes the proof. □

THEOREM 5.4. *Under the assumptions of Theorem 5.3 let $\varphi_{(n)}$ be an SUD test of order $\lambda_n$ with $\liminf_{n\to\infty} \lambda_n/n > 0$ and the condition (5.5) be replaced by*

$$(5.8) \qquad \forall \varepsilon > 0 : \liminf_{n\to\infty} \inf_{\varepsilon \leq t \leq K}(t - \hat{F}_n(\rho(t))) > 0 \quad [P_\vartheta]$$

*for some $K \in [0,1]$ fulfilling $K > L = \limsup_{n\to\infty} \lambda_n/n$ in the case of $L < 1$ and $K = 1$ otherwise. Supposing that $\lim_{x \downarrow 0} q(x)$ exists, we have*

$$(5.9) \qquad \lim_{n\to\infty} \mathrm{FDR}_\vartheta(\varphi_{(n)}) = \zeta \lim_{x \downarrow 0} q(x) = \zeta q(0) = \zeta \overline{q}(0) \leq \zeta \alpha.$$

PROOF. Again, to avoid triviality, we assume $I_0(\vartheta) \neq \varnothing$. Equation (5.9) can be shown by utilizing the notation introduced in the proof of Theorem 4.3 and the decomposition

$$\mathrm{FDR}_\vartheta(\varphi_{(n)}) = \zeta_n \sum_{j=1}^{\lambda_n} q(j/n) P_\vartheta(R_n = j | p_{i_0} \leq \alpha_{j:n})$$
$$+ \zeta_n \sum_{j=\lambda_n+1}^{n} q(j/n) P_\vartheta(R_n = j | p_{i_0} \leq \alpha_{j:n})$$
$$= M_n + m_n \qquad \text{(say)}.$$

In view of Theorem 4.3 and the structure of an SUD test, we obtain by applying (4.2) that

$$M_n = \zeta_n \mathbb{E}_{\vartheta^{i_0}}[q(R_n/n)\mathbf{1}_{\{R_n/n \leq \lambda_n/n\}}],$$
$$m_n \leq \zeta_n \mathbb{E}_{\vartheta^{i_0}}[\overline{q}(R_n/n)\mathbf{1}_{\{R_n/n > \lambda_n/n\}}].$$

From (5.8) it follows that $P_\vartheta$-almost surely $\hat{F}_n(\rho(\lambda_n/n)) < \lambda_n/n \leq K$ and consequently $R_n/n < \lambda_n/n \leq K$ holds true for eventually all $n \in \mathbb{N}$. Therefore, again due to (5.8), in analogy to the proof of Theorem 5.3 we conclude that $\lim_{n\to\infty} R_n/n = 0$ $P_{\vartheta^{i_0}}$-almost surely, which finally entails $\lim_{n\to\infty} \mathbf{1}_{\{R_n/n > \lambda_n/n\}} = 0$ $P_{\vartheta^{i_0}}$-almost surely. Together with the boundedness of $\overline{q}$ this entails that $\lim_{n\to\infty} m_n = 0$. Moreover, exploiting the continuity of $q$ at $x = 0$ we see that $\lim_{n\to\infty} M_n = \zeta q(0) = \zeta \overline{q}(0)$, which altogether yields the desired result. □

REMARK 5.1. One cannot expect to obtain exact values for the limiting FDR under the quite general assumptions of Theorem 5.2 if $r^* = 0$. To see this, consider the case $\zeta_n \equiv 1$ in which the FDR is equal to the familywise error rate (FWER). For $\zeta_n \equiv 1$ it was shown in [12] that the FWER is equal to $\alpha$ for any $n \in \mathbb{N}$ in the case of a linear SU procedure, while it tends to $1 - \exp(-\alpha) < \alpha$ for a linear SD procedure. We therefore have to know more about the structure of the underlying procedure in order to compute the limiting FDR in case of $r^* = 0$. The limiting behavior for procedures based on $f_\alpha$ (or its modifications) satisfying the assumptions of Theorem 5.4 is in accordance with the linear SU procedure and should be expected, since the difference of the critical values $\alpha_{i:n} - i\alpha/n$ tends to zero for $i \in o(n)$. Therefore, the local behavior around zero should not differ much for large $n$.

COROLLARY 5.1 (Examples 3.1–3.3 continued). *Assume the distributional assumptions (D3), (I1) and (I2) hold. Then the SUD procedure based on $f_\alpha$ with parameter $\lambda \in [0,1)$ and the SU procedures based on $f_{\alpha,\kappa}^{(i)}$, $i = 1, 2$, as well as the truncated SU procedure asymptotically control the FDR at level $\alpha$. More precisely, if condition (5.3) is fulfilled, that is, $\lim_{n\to\infty} \zeta_n = \zeta \in [0,1]$, then:*

(i) *For the SUD procedure the upper bound $\alpha$ for the limiting FDR is sharp for $\zeta \in [\alpha, 1]$.*

(ii) *For the SU procedures based on $f_{\alpha,\kappa}^{(i)}$, $i = 1, 2$, the upper bound $\alpha$ for the limiting FDR is sharp for $\zeta \geq \zeta^*(\kappa) = \alpha/(\kappa(1-\alpha)+\alpha)$. In the case of $\zeta < \zeta^*(\kappa)$, an upper bound for the asymptotic FDR is given by $\zeta \tilde{t}_\zeta/(1 - \zeta + \zeta \tilde{t}_\zeta)$, where $\tilde{t}_\zeta$ denotes the unique solution of the equation $F_\infty(t|\zeta) = h_i(t), i = 1, 2$, on $(0, t_\zeta)$. For finite $n$, the upper bound given in (4.9) is sharp.*

(iii) *For the truncated SU procedure the upper bound $\alpha$ for the limiting FDR is sharp for $\zeta \geq \zeta^*(\kappa)$. In the case of $\zeta < \zeta^*(\kappa)$, an upper bound for the asymptotic FDR is given by $\zeta\kappa/(1 - \zeta + \zeta\kappa)$.*

PROOF. First of all, as mentioned before, a step-up-down test has the structural properties (T1), (T2) and (A1). Moreover, assumptions (D3), (I1) and (I2) imply the crucial monotonicity properties (D1) and (D2) for a step-up-down test. Hence, in order to apply Theorem 5.1, it remains to check the validity of condition (5.1). To this end, for notational convenience and without loss of generality, we work under condition (5.3). We make use of (2.1), that is, that the ecdf $\hat{F}_n$ converges $P_{I'_0}$-almost surely to its limit $F_\infty(\cdot|\zeta)$ uniformly in $t \in [0,1]$. Since under (T1) and (T2) we have the identity $\hat{F}_n(\rho(R_n/n)) = R_n/n$ for all $n \in \mathbb{N}$, (2.1) leads to $\lim_{n\to\infty}(F_\infty(\rho(R_n/n)|\zeta) - R_n/n) = 0$ $P_{I'_0}$-almost surely. From this we conclude that ($P_{I'_0}$-almost surely) the only possible accumulation points of the sequence $(R_n/n)_{n\in\mathbb{N}}$ consist of the solutions of the equation $F_\infty(\rho(t)|\zeta) = t$ in $t \in [0,1]$. If, as in Examples 3.2 and 3.3, this solution is unique, then the sequence $(R_n/n)_{n\in\mathbb{N}}$ necessarily converges to this solution $r^* = r^*(\zeta) \in [0,1]$ (say) $P_{I'_0}$- almost surely.

If, however, as in Example 3.1, the equation $F_\infty(\rho(t)|\zeta) = t$ has the solution $t = 1$ and exactly one further (smaller) solution in $[0,1)$, we have to exclude $t=1$ as a possible accumulation point of $(R_n/n)_{n\in\mathbb{N}}$ in the latter case in order to prove the $P_{I_0'}$-almost sure convergence of $(R_n/n)_{n\in\mathbb{N}}$ to the smallest solution $r^* = r^*(\zeta)$ (say) of the aforementioned equation. To this end, we only consider values of $\zeta$ leading to the two distinct solutions $r^*(\zeta)$ and 1. For critical value functions $\rho$ with $\rho(t) \leq f_\alpha^{-1}(t)$ for all $t \in [0,1]$ it is then evident that $F_\infty(\rho(t)|\zeta) < t$ for all $t \in (r^*, 1)$. Moreover, notice that, by definition of $\lambda_n(\lambda)$, we have the inequalities $\lambda_n(\lambda) - 1 \leq nr(\lambda) \leq \lambda_n(\lambda)$ for all $n \in \mathbb{N}$. Now, if $\lambda > \rho(r^*)$, this, together with condition (2.1), yields that $P_{I_0'}$-almost surely $\hat{F}_n(\rho(\lambda_n(\lambda)/n)) < \lambda_n(\lambda)/n$ and consequently $R_n < \lambda_n(\lambda)$ holds true for eventually all $n \in \mathbb{N}$. This entails $\limsup_{n\to\infty} R_n/n \leq \limsup_{n\to\infty} \lambda_n(\lambda)/n = r(\lambda) < 1$ $P_{I_0'}$-almost surely, which is just what we wanted to show.

Finally, if $\lambda \leq \rho(r^*)$, we may choose a $\lambda' > \rho(r^*)$ and compare the number of rejections of the corresponding SUD procedures. Since this number is non-decreasing in $\lambda \in [0,1]$, we eventually arrive at $\limsup_{n\to\infty} R_n/n \leq r(\lambda') < 1$ $P_{I_0'}$-almost surely.

Since for all procedures under investigation it holds $\rho(t) \leq f_\alpha^{-1}(t)$ for all $t \in [0,1]$, we conclude that $r^* = r^*(\zeta) \leq f_\alpha(t_\zeta)$. Hence, Theorem 5.1 applies. As a consequence, the example procedures asymptotically control the FDR.

In the case of the SUD procedure in (i), we use $\rho = f_\alpha^{-1}$ and obtain $r^* = f_\alpha(t_\zeta)$. Hence, the upper bound $\alpha$ for the asymptotic FDR is sharp in (i) under Dirac-uniform configurations. The sharpness of the upper bound $\alpha$ for the asymptotic FDR in (ii) and (iii) is due to the fact that under Dirac-uniform configurations with $\zeta \geq \zeta^*(\kappa)$ we obtain $r^* = f_\alpha(t_\zeta)$.

Finally, the sharpness of the upper bounds for the finite $n$ FDR in (ii) is a consequence of (A2), which is fulfilled for $f_{\alpha,\kappa}^{(i)}$, $i=1,2$. Sharpness here means that the upper bound given in (4.9) is exactly attained for finite $n$ under Dirac-uniform configurations. $\square$

The latter corollary means, in other words, that procedures based on $f_\alpha$ fulfilling the assumptions of Theorem 5.2 asymptotically exhaust the whole FDR level $\alpha$ under Dirac-uniform configurations. Moreover, the rejection curve $f_\alpha$ cannot be improved in the sense of the following theorem, which is another consequence of Theorem 5.1. In order to formalize this, let $\alpha \in (0,1)$, $\lambda \in [0,1]$ and $\mathcal{M}_\lambda$ denote the set of rejection curves $r$ with the property that for all $I_0 \subseteq \mathbb{N}$ with $\lim_{n\to\infty} \zeta_n = \zeta$ for some $\zeta \in [0,1]$ it holds

$$(5.10) \quad \limsup_{n\to\infty} \mathrm{FDR}_{I_0}(\varphi_{n,\lambda_n}^{\mathrm{SUD}(r)}) \leq \limsup_{n\to\infty} \sup_{\vartheta \in \Theta} \mathrm{FDR}_\vartheta(\varphi_{n,\lambda_n}^{\mathrm{SUD}(r)}) \leq \alpha,$$

where $\varphi_{n,\lambda_n}^{\mathrm{SUD}(r)}$ is the step-up-down procedure of order $\lambda_n = \lambda_n(\lambda)$ based on $r$. It should be noted that the first inequality in (5.10) is not very restrictive since many statistical models satisfy the "model continuity assumptions" (SA) formulated in [11], due to which, at least for SUD procedures such as $\varphi_{n,\lambda_n}^{\mathrm{SUD}(r)}$, the corresponding FDR values $\mathrm{FDR}_{I_0}(\varphi_{n,\lambda_n}^{\mathrm{SUD}(r)})$ can be approximated arbitrarily closely by the values $\mathrm{FDR}_\vartheta(\varphi_{n,\lambda_n}^{\mathrm{SUD}(r)})$ for some suitably chosen $\vartheta \in \Theta$, $n \in \mathbb{N}$.

In terms of power it is immediately clear that, whenever $r_1, r_2 \in \mathcal{M}_\lambda$ with $r_1 \leq r_2$, then $\varphi_{n,\lambda_n}^{\mathrm{SUD}(r_1)} \geq \varphi_{n,\lambda_n}^{\mathrm{SUD}(r_2)}$. Therefore, a smaller rejection curve typically leads to a more powerful test procedure in the sense that more (never less) false hypotheses can be rejected. Here we define the power of a test $\varphi_{(n)}$ by $\bar{\beta}_\vartheta(\varphi_{(n)}) = \mathbb{E}_\vartheta[(R_n - V_n)/(n_1 \vee 1)]$.

THEOREM 5.5 (Asymptotic optimality of $f_\alpha$).

(i) Let $\lambda \in [0,1]$ and $r \in \mathcal{M}_\lambda$. Then

(5.11) $$\forall t \in [0,\lambda] : r(t) \geq f_\alpha(t).$$

If $\lambda < 1$, then it holds for any $\tau \in (\lambda, 1]$ that

(5.12) $$\forall t \in (\lambda, \tau] : r(t) \leq f_\alpha(t) \Rightarrow \forall t \in (\lambda, \tau] : r(t) = f_\alpha(t).$$

(ii) If $\lambda < 1$ and $r \in \mathcal{M}_\lambda$ is such that, for every $\zeta \in (\alpha, 1)$, the equation $F_\infty(\rho(t)|\zeta) = 1 - \zeta + \zeta\rho(t) = t$ has at most one solution in (0, 1), it even holds $r(t) \geq f_\alpha(t)$ for all $t \in [0,1]$.

(iii) If $\lambda = 1$ and assuming (D3), (I1) and (I2), it holds that

$$\inf_{r \in \mathcal{M}_1} r = f_\alpha.$$

Moreover, for any $\vartheta \in \Theta_\kappa = \{\vartheta \in \Theta : \liminf_{n \to \infty} \zeta_n(\vartheta) > \alpha/(\kappa(1-\alpha)+\alpha)\}$, $\kappa \in (0,1)$, the power of any $\tilde{f}_\alpha \in \mathcal{M}_1$ with $\tilde{f}_\alpha(t) = f_\alpha(t)$ for all $t \in [0,\kappa]$ is asymptotically not smaller than the power of any other $r \in \mathcal{M}_1$, that is,

(5.13) $$\liminf_{n \to \infty} [\bar{\beta}_\vartheta(\varphi_{n,n}^{\mathrm{SUD}(\tilde{f}_\alpha)}) - \bar{\beta}_\vartheta(\varphi_{n,n}^{\mathrm{SUD}(r)})] \geq 0 \qquad \textit{for all } \vartheta \in \Theta_\kappa.$$

PROOF. In order to prove part (i), assume that for an arbitrary chosen rejection curve $r \in \mathcal{M}_\lambda$ it holds $r(t^*) < f_\alpha(t^*)$ for some $t^* \in (0,\lambda)$. Consider now a Dirac-uniform configuration $P_{I_0}$ with $\lim_{n \to \infty} \zeta_n = \zeta$ and $\zeta \in (\alpha, 1)$ chosen such that $r(t^*) < F_\infty(t^*|\zeta) < f_\alpha(t^*)$. Then it is obvious that property (5.10) is violated, because (with self-explaining notation) it follows $P_{I_0}$-almost surely

$$\liminf_{n \to \infty} R_n^{(r)}/n \geq F_\infty(t^*|\zeta) > F_\infty(t_\zeta|\zeta) = f_\alpha(t_\zeta)$$

and consequently

$$\liminf_{n \to \infty} \mathrm{FDR}_{I_0}(\varphi_{n,\lambda_n}^{\mathrm{SUD}(r)}) \geq \zeta t^*/(1 - \zeta + \zeta t^*) > \zeta t_\zeta/(1 - \zeta + \zeta t_\zeta) = \alpha,$$

due to the fact that the function $x \mapsto \zeta x/(1 - \zeta + \zeta x)$ is strictly increasing in $x \in (0,1)$ and $t^* > t_\zeta$. Hence, for all $t \in (0, \lambda)$ we have $r(t) \geq f_\alpha(t)$, from which the assertion follows.

Now assume that we have $r(t) \leq f_\alpha(t)$ for all $t \in (\lambda, \tau]$ and $r(t^*) < f_\alpha(t^*)$ for some $t^* \in (\lambda, \tau)$. Consider now the Dirac-uniform asymptotic model $\mathrm{DU}_\infty(\zeta^*)$ with $\zeta^* \in (\alpha, 1)$ chosen such that $f_\alpha(\lambda) < F_\infty(\lambda|\zeta^*)$, $F_\infty(t^*|\zeta^*) < f_\alpha(t^*)$ and $\inf_{\lambda \leq t \leq t^*}(F_\infty(t|\zeta^*) - r(t)) > 0$, which is possible due to the left-continuity of the rejection curve $r$. Then the argumentation is the same as before. Part (ii) and the first assertion of part (iii) can be proven similarly.

For the proof of (5.13), we assume (in order to avoid triviality) $n_1(n) > 0$ for all $n \in \mathbb{N}$, define $S_n = R_n - V_n$ and denote the set of all $\tilde{f}_\alpha \in \mathcal{M}_1$ with $\tilde{f}_\alpha(t) = f_\alpha(t)$ for all $t \in [0, \kappa]$ by $\mathcal{S}_\kappa$. Then we have (with self-explaining notation as before) the inequality

$$\forall n \in \mathbb{N} : \forall \tilde{f}_\alpha \in \mathcal{S}_\kappa : \forall r \in \mathcal{M}_1 : \left(\frac{S_n(\tilde{f}_\alpha)}{n_1} - \frac{S_n(r)}{n_1}\right) \mathbf{1}_{\{t_n^*(r) \leq \kappa\}} \geq 0,$$

which holds true due to (5.11) and the fact that $S_n$ is nondecreasing in $t_n^*$. Now, for fixed $\vartheta \in \Theta_\kappa$, we utilize the chain of inequalities

$$t_n^*(r|P_\vartheta) \leq t_n^*(r|\mathrm{DU}_n(\zeta_n(\vartheta))) \leq t_n^*(\tilde{f}_\alpha|\mathrm{DU}_n(\zeta_n(\vartheta))) < \kappa$$

which holds $P_\vartheta$-almost surely for eventually all $n \in \mathbb{N}$, leading to $\limsup_{n\to\infty} t_n^* (r|P_\vartheta) < \kappa$ and consequently to $\mathbf{1}_{\{t_n^*(r) \leq \kappa\}} \to 1\ [P_\vartheta]$ for all $\vartheta \in \Theta_\kappa$. Therefore, we obtain $P_\vartheta$-almost surely

$$(5.14) \quad \liminf_{n \to \infty} \left(\frac{S_n(\tilde{f}_\alpha)}{n_1} - \frac{S_n(r)}{n_1}\right) \geq 0 \quad \text{for all } \vartheta \in \Theta_\kappa, \tilde{f}_\alpha \in \mathcal{S}_\kappa, r \in \mathcal{M}_1.$$

Taking expectation in (5.14) and utilizing Fatou's lemma, we finally arrive at assertion (5.13). □

Theorem 5.5 shows that in the class of SU procedures with rejection curve $r \in \mathcal{M}_1$ we always have $r \geq f_\alpha$. In the class of truncated SU procedures with parameter $\kappa \in (0,1)$, the truncated procedure based on $f_\alpha$ is the best choice. More generally, if we restrict attention to the subspace $\Theta_\kappa \subset \Theta$ described in (iii) of Theorem 5.5, $f_\alpha$ is the asymptotically uniformly best choice on $[0,\kappa]$ for a step-up procedure. For SUD procedures with parameter $\lambda < 1$, $f_\alpha$ leads to the asymptotically uniformly best choice of critical values on the step-up part, see (5.11). On the step-down part of a SUD procedure, $f_\alpha$ cannot be uniformly improved by some $r \in \mathcal{M}_\lambda$ whatever $r$ does on the step-up part; see (5.12) with $\tau = 1$. For arbitrary $\tau \in (\lambda, 1]$, assertion (5.12) states that a rejection curve $r \in \mathcal{M}_\lambda$ cannot be first smaller and then larger than $f_\alpha$ on the interval $(\lambda, 1]$. It seems possible that $\mathcal{M}_\lambda$ contains an $r$ which is first larger and then smaller on the step-down part. But this would imply that the SUD procedure based on $r$ is asymptotically less powerful than the SUD procedure based on $f_\alpha$ on some $\Theta_\kappa$. If we restrict attention to rejection curves $r \in \mathcal{M}_\lambda$ described in (ii) of Theorem 5.5, then $f_\alpha$ is the best choice. These considerations may justify calling $f_\alpha$ the asymptotically optimal rejection curve (AORC).

**6. Concluding remarks.** In view of the asymptotic optimality results developed in Section 5 concerning procedures based on $f_\alpha$ or its modifications, it is natural to ask how large $n$ has to become in order to achieve a reasonable behavior of the FDR of the proposed procedures. As already mentioned in Example 3.1, the asymptotic exhaustion of the whole FDR level has to be traded off with a slightly liberal behavior of the procedures based on $f_\alpha$ in the finite case. In order to illustrate this effect, we consider the SU procedures based on $f^{(i)}_{\alpha,\kappa}$, $i = 1, 2$, where the upper bound given in (4.9) is sharp in the $\mathrm{DU}_n(\zeta_n)$-model. Due to the pointwise order of these two rejection curves (see Figure 1) it is clear that an SU procedure based on $f^{(2)}_{\alpha,\kappa}$ is more liberal in the $\mathrm{DU}_n(\zeta_n)$-model. We therefore present results for this procedure. Figure 2 depicts the behavior of this procedure under Dirac-uniform configurations with a varying number of true hypotheses for three different values of $n$. For $n = 100$, there is a notable violation of the FDR level $\alpha = 5\%$ for $12 \leq n_0 \leq 35$. The largest FDR under Dirac-uniform is attained for $n_0 = 16$ with numerical value 0.05801. For the two larger values of $n$, the actual level does not exceed $\alpha$ by much. Computation of the FDR for a SU(D) procedure in the case of an underlying Dirac-uniform configuration can be done by utilizing formulas for the joint distribution function of order statistics, see [21], pages 366–367, and [10].

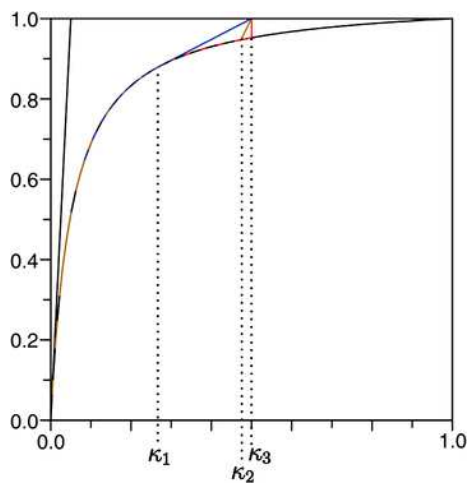

FIG. 1. From left to right: Simes line, $f^{(1)}_{\alpha,\kappa_1}$, $f^{(2)}_{\alpha,\kappa_2}$, the truncated version based on $f_\alpha$ with $\kappa_3 = 1/2$ and $f_\alpha$. The $\kappa_i$'s are chosen such that $f^{(i)}_{\lambda,\kappa_i}(1/2) = 1$ for $i = 1, 2$.

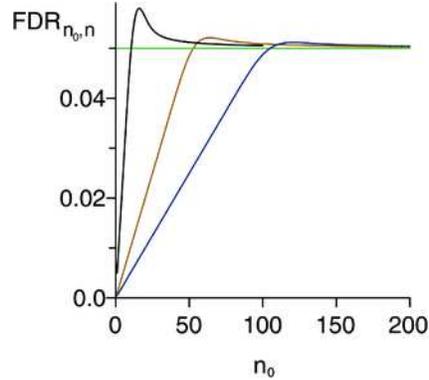

Fig. 2. *Actual FDR of the SU procedure based on $f^{(2)}_{0.05,\kappa_2}$ for $n = 100, 500, 1000$ depending on $n_0$ under Dirac-uniform configurations. The curves can be distinguished by noticing that the maximum FDR becomes smaller for increasing $n$.*

We will give one brief suggestion for a modification of $f_\alpha$ in the finite case. However, this shall not be emphasized too much, because on the one hand, the AORC is designed for the asymptotic case and on the other hand, we have to keep in mind that the FDR values under Dirac-uniform reflect an unrealistic worst case scenario. For realistic alternatives, we get much smaller realized FDRs so that the original AORC may safely be used in the finite case, for example, $n \geq 500$.

We only mention one possibility to obtain a valid set of critical values for an SU or SUD procedure guaranteeing strict FDR control, that is, we adjust the critical values given in (3.1) in an appropriate way. For example, we can try to find a suitable $\beta_n > 0$ such that the choice

$$(6.1) \qquad \alpha_{i:n} = \frac{i\alpha}{n + \beta_n - i(1-\alpha)}, \qquad i = 1,\ldots,n,$$

yields an SU procedure (or SUD procedure) controlling the FDR at level $\alpha$. The critical values (6.1) correspond to the rejection curve

$$\tilde{f}_\alpha(t) = \left(1 + \frac{\beta_n}{n}\right) f_\alpha(t), \qquad t \in [0, \alpha/(\alpha + \beta_n/n)].$$

It is remarkable that herewith a direct connection to the considerations in [3] can be drawn. In the remark to Definition 7 in [3], the authors propose an (adaptive) SD procedure with critical values given by (6.1) and the universal adjustment constant $\beta_n \equiv 1.0$. FDR control for this SD procedure is proved in [15] in case of independent $p$-values. For an SU procedure, the adjustment constant has to be larger. For example, for $\alpha = 0.05$, an SU procedure with $n = 100$ and the choice $\beta_{100} = 1.76$ leads to strict FDR control.

A systematic comparison of other procedures controlling the FDR strictly or asymptotically with procedures based on the AORC goes beyond the scope of this paper and is directed to future research. Typically, each approach has assets and drawbacks. At least in the class of SUD procedures (including SD and SU procedures) based on fixed rejection curves the new procedures based on the AORC are powerful alternatives. One can construct situations, especially if the proportion $\zeta$ of true hypotheses is small, where the new procedures reject many more hypotheses than the LSU procedure. On the other hand, if $\zeta$ is large, there may be only a few more or, in rare cases, fewer rejections. In other words, if $\zeta$ is close to one, it is hard to beat the LSU procedure by a considerable amount. Moreover, it should not be concealed that the LSU procedure has the advantage that it applies for certain situations with positive dependent $p$-values, whereas procedures based on the AORC may fail to control the FDR in such situations. For example, this is the case for $p$-values based on normally distributed test statistics with positive correlation as it appears in multiple comparisons with a control. Whether the AORC works for pairwise comparisons or multiple testing of correlation coefficients $\rho_{ij}$ ($H_{ij}: \rho_{ij} = 0$ versus $K_{ij}: \rho_{ij} \neq 0$, $1 \leq i < j \leq k$) is currently under investigation.

Finally, one may think about a more flexible concept of FDR control depending on the proportion $\zeta_n = \zeta_n(\vartheta)$ of true hypotheses which allows, for example, a larger FDR for larger values of $\zeta_n$ and a smaller FDR for smaller values of $\zeta_n$, or vice versa. Therefore, one may choose a suitable FDR controlling function $g:[0,1] \longrightarrow [0,1]$ and require FDR control at this level function $g$, that is,

$$\forall \vartheta \in \Theta : \mathrm{FDR}_\vartheta(\varphi_{(n)}) \leq g(\zeta_n(\vartheta)).$$

This means that for any $\zeta_n \in [0,1]$ the FDR is controlled at level $g(\zeta_n)$ if the proportion of true hypotheses is $\zeta_n$. Obviously, for the LSU procedure at FDR level $\alpha$ we can choose $g(\zeta) = \zeta\alpha$. For an SUD procedure related to the AORC with critical values defined in (6.1) and suitable $\beta_n$ we can choose $g(\zeta) = \min\{\alpha, \zeta\}$. Asymptotic FDR control at level function $g$ now means that

$$\limsup_{n \to \infty} \sup_{\vartheta \in \Theta} (\mathrm{FDR}_\vartheta(\varphi_{(n)}) - g(\zeta_n(\vartheta))) \leq 0.$$

For example, in the case of the truncated SU procedure, we get from Corollary 5.1 that we can choose $g(\zeta) = \alpha$ for $\zeta \in [\zeta^*(\kappa), 1]$ and $g(\zeta) = \zeta\kappa/(1 - \zeta + \zeta\kappa)$ for $\zeta \in [0, \zeta^*(\kappa))$.

**Acknowledgments.** The authors are grateful to the referees, an Associate Editor and the editor B. W. Silverman for their constructive and valuable comments and suggestions.

H. FINNER
T. DICKHAUS
GERMAN DIABETES CENTER
INSTITUTE OF BIOMETRICS AND EPIDEMIOLOGY
DÜSSELDORF
GERMANY
E-MAIL: finner@ddz.uni-duesseldorf.de
  dickhaus@ddz.uni-duesseldorf.de

M. ROTERS
BIOMETRICS DEPARTMENT
OMNICARE CLINICAL RESEARCH
KÖLN
GERMANY
E-MAIL: Markus.Roters@omnicarecr.com